\documentstyle[12pt, twoside]{article}
\input amssym.def
\input amssym.tex 

\newenvironment{demo}[1]%
{\vskip-\lastskip\medskip
  \noindent
  {\em #1.}\enspace
  }%
{\qed\par\medskip
  }

\newcommand{\qed}{
  \strut\hfill
  \mbox{$\Box$}
  }

{\vskip-\lastskip
  \begin{trivlist}
    \item[]
      {\bf Axiom #1 }
      }%
{\end{trivlist}
  }

\newtheorem{theorem}{Theorem}[section]

\newtheorem{remark}{Remark}[section]

\newtheorem{proposition}{Proposition}[section]

\newcommand{\bb}{
  {\bf b}
  }

\newcommand{\bc}{
  {\bf c}
  }

\newcommand{\bFmn}{
  { {\cal F}_{\pm} }^{\otimes (M+N)}
  }

\newcommand{\bFm}{
  { {\cal F}_{\pm} }^{\otimes M}
  }

\newcommand{\bFn}{
  { {\cal F}_{\pm} }^{\otimes N}
  }

\newcommand{\C}{
  \Bbb C
  }

\newcommand{\D}{
  \cal D
  }

\newcommand{\F}{
  \cal F
  }

\newcommand{\FBn}{
  {\cal F}_{-}^{\otimes N}
  }

\newcommand{\FFn}{
  {\cal F}_{+}^{\otimes N}
  }

\newcommand{\g}{
  \frak g
  }

\newcommand{\gl}{
  {gl}
  }

\newcommand{\hD}{
  \widehat{\cal D}
  }

\newcommand{\hg}{
  \widehat{\frak g}
  }

\newcommand{\hgl}{
  \widehat{gl}
  }

\newcommand{\hL}{
  \widehat{\Lambda}
  }

\newcommand{\Lp}{
  {L_N (\Lambda_{+} (\lambda) ) }
  }

\newcommand{\N}{
  (n)z^{-n}
  }

\newcommand{\NN}{
  (n)z^{-n-1}
  }

\newcommand{\Of}{
  {\cal O}_{f}
  }

\newcommand{\Ofm}{
  {\cal O}_{f}^M
  }

\newcommand{\Ofn}{
  {\cal O}_{f}^N
  }

\newcommand{\OpN}{
  {\cal O}^N_{+}
  }

\newcommand{\Op}{
  {\cal O}_{+}
  }

\newcommand{\Om}{
  {\cal O}_{-}
  }

\newcommand{\Opm}{
  {\cal O}_{\pm}
  }

\newcommand{\OWp}{
  {\cal O}^w
  }

\newcommand{\pairpm}{
  ( GL(N), \hgl |_{c= \pm N} )
  }

\newcommand{\SUM}{
  \sum_{n\in \Bbb Z}
  }

\newcommand{\Tr}{
  \mbox{Tr}
  }

\newcommand{\vac}{
  |0 \rangle
  }

\newcommand{\W}{
  { \cal W }_{1+\infty}
  }

\newcommand{\Wth}{
  {\cal W}_3
  }

\newcommand{\WGN}{
  {\cal{W} } (gl_N)
  }

\newcommand{\Z}{
  \Bbb Z
  }

\begin{document}
\title{
  Dual pairs and tensor categories of modules 
  over Lie algebras $\widehat{gl}_{\infty}$ and $\W$
  }
\author{
  %
  Weiqiang Wang\\
\\{\small Department of Mathematics, 
Yale University, New Haven, CT 06520}}
\date{}

\maketitle
\begin{abstract}
  We introduce a tensor category $\Op$ (resp. $\Om$) of
certain modules of $\widehat{gl}_{\infty}$ with non-negative
(resp. non-positive) integral central charges with the usual
tensor product.
We also introduce a tensor category $(\Of, \bigodot)$
consisting of certain modules over $GL(N)$ for all $N$.
We show that the tensor categories $\Opm$ and $\Of$ are 
semisimple abelian and all equivalent to each other.
We give a formula to decompose a tensor product of
two modules in each of these categories. We also introduce
a tensor category $\OWp$ of certain modules over $\W$
with non-negative integral central charges. We show
that $\OWp$ is semisimple abelian and give an explicit
formula to decompose a tensor product of
two modules in $\OWp$.
\end{abstract}
\setcounter{section}{-1}
\section{Introduction}

A general important problem in representation
theory is to study the decomposition
of a tensor product of two representations
of a Lie algebra $\g$.
For example, due to Weyl's complete reducibility
theorem, a tensor product of two finite dimensional
representations of a complex simple
Lie algebra $\g$ is decomposed into a direct
sum of irreducibles.

However when one considers a similar problem
in the representation theory of infinite dimensional
Lie algebras, the story is usually much more complicated.
For example, if one decomposes
the usual tensor product of an ireducible integrable 
representation with a positive integral central charge
with another irreducible integrable representation
with another positive integral central charge of
an affine Kac-Moody algebra $\hg$ \cite{K}, one
gets an infinite direct sum of irreducibles with infinite
multiplicities. If one considers a similar question
for the negative central charge cases, one usually
cannot expect to get any reasonable decomposition due to
the failure of complete reducibility.

Lie algebra $\widehat{gl}_{\infty}$,
which throughout this paper will be 
denoted by $\hgl$ for shortness,
plays an important role in the
representation theory of infinite dimensional Lie
algebras, with remarkable relations with solitons,
KP equations, affine Kac-Moody algebras 
and the $\W$ algebra etc
(cf. eg. \cite{DJKM}, \cite{KRa}, \cite{KR1}).
In this paper, we introduce two tensor categories
$\Op$ and $\Om$ of representations of $\hgl$.
The irreducible objects in  category
$\Op$ (resp. $\Om$) are certain irreducible
highest weight representations with
non-negative (resp. non-positive)
integral central charges. Note that central
charges in each category $\Opm$ are not fixed.

We show that a tensor product of two
irreducible modules in category $\Op$ (resp. $\Om$)
decomposes into an (infinite)
direct sum of irreducibles with
FINITE multiplicities. We identify the
irreducibles which appear in such a decomposition 
and give the multiplicities of these irreducibles
in terms of certain
multiplicities of irreducible representations
appearing in the decomposition of certain irreducible
finite dimensional representations of a general linear 
Lie group when
restricted to some semisimple block-diagonal subgroup.
This is indeed a reciprocity law
of branching rules due to the
appearance of the so-called ``seesaw'' pairs \cite{Ku, H2}.

We then show that the two categories $\Opm$ 
are semisimple abelian tensor categories
and are equivalent to each other.
In the category $\Op$ case,
an important ingredient is a remarkable
duality theorem first
proved by I.~Frenkel \cite{F1, FKRW} among representations
of $\hgl$ with positive integral
central charge $N$ and those of $GL(N)$.
In the category $\Om$ case, a similar
duality theorem from \cite{KR2} is used.
We also define $\Ofn$ to be the category of all 
(possiblly infinite dimensional) representations
whose Jordan-Holder series consisting of finite dimensional
irreducible $GL(N)$-modules. Denote by $\Of$ the category which is
a direct sum of categories $\Ofn$ for all $N \geq 0$.
We introduce a tensor product $\bigodot$ on the category $\Of$.
and then establish an equivalence of 
tensor categories between $\Of$ and $\Opm$.

In \cite{KR1}, Kac and Radul systematically developed 
the theory of quasifinite 
representations of $\hD$ by using
connections between $\hD$ and
$\hgl$ (or rather some variations of $\hgl$) in an ingenious way.
Here $\hD$ is the universal central extension of
Lie algebra of differential operators on the circle, which
is often referred to as $\W$ in literature.
In this paper we define a semisimple abelian
tensor category $\OWp$ of certain
highest weight $\hD$ modules
with non-negative integral central charge. 
We prove that $\OWp$ is a semisimple abelian
tensor category. Indeed
this category consists of exactly all the
representations of the vertex algebra $\W$ 
with non-negative integral central charges \cite{FKRW}.
Based on our tensor product decomposition theorem in category $\Op$
and relations between $\hgl$ and $\hD$, an
explicit tensor product decomposition 
in the category $\OWp$ is presented.
We remark that a similarly defined category of 
$\hD$-modules with non-positive
integral central charge will
not be a semisimple tensor category.

Let us explain in more detail. In the case of positive
(resp. negative) integral central charges, 
we will need $N$ pairs of
fermions $b^p(z), c^p (z)$
(resp. bosonic ghosts $\beta^p(z), \gamma^p (z) $), 
$p =1, \dots, N. $
We will use bold letters when we want to treat both cases 
at the same time. For instance,
we use $\bb (z)$ (resp. $\bc (z)$) to represent
either $b(z)$ or $\beta(z)$ (resp. $c(z)$ or $\gamma (z)$).
Denote by $\F_{+}$ (resp. $\F_{-}$) the
Fock space of a pair of fields $b(z), c(z)$ 
(resp. $\beta (z), \gamma (z)$).
Now take $N$ pairs of $\bb\bc$ fields and
consider the tensor product $\bFn$, namely the Fock space of $N$ pairs
of $\bb\bc$ fields. Lie algebras $\hgl$ and
$gl(N)$ acting on $\bFn$ naturally commute with each other
and form a dual pair.
The action of $gl(N)$ can be lifted to an action of $GL(N)$.
Then with respect to the commuting
actions of the dual pair $(GL(N), \hgl)$, 
$\bFn$ decomposes into a direct sum of
irreducible isotypic subspaces 
(cf. \cite{F1}, \cite{FKRW}, \cite{KR2}, 
see section \ref{sect_decomp} for notations):
$$ \bFn = \bigoplus_{\lambda \in P^N_{+} } 
V_{\lambda} \bigotimes L(\Lambda_{\pm} (\lambda)), $$
where $V_{\lambda}$ is the irreducible highest weight
$GL(N)$-module of highest weight $\lambda$, and 
$L_{\pm N} (\Lambda_{\pm} (\lambda) )$ is the irreducible
$\hgl$-module with highest weight $\Lambda_{\pm} (\lambda) $
and central charge $\pm N$.
We often write $(GL(N), \hgl)$ as 
$\left( GL(N), \hgl |_{c= \pm N} \right)$ 
to emphasize that the action of $\hgl$ on $\bFn$ has
central charge $\pm N$. Keeping the obvious isomorphism
between $\bFm \otimes \bFn$ and $\bFmn $ in mind, we
see that the following two dual pairs
act on the same Fock space $\bFmn$:

\begin{eqnarray*}
 \left\{ \begin{array}{cc}
  \Bigl( GL(M+N), 
    & 
    \hgl |_{c=\pm (M+N) } \Bigl) \nonumber \\
 \uparrow                          & \downarrow \nonumber\\
 \Bigl( GL(M)\times GL(N), 
     & \hgl |_{c=\pm (M) } \bigoplus \hgl |_{c=\pm (N) } \Bigl)
                       \end{array} \right. 
\nonumber
\end{eqnarray*}
where inclusions of Lie groups/algebras are shown by the arrows. 
Thus these
two dual pairs form a seesaw pair \cite{Ku, H2}. A consequence
of a seesaw pair will be a duality among the branching rules
corresponding to the two inclusions of Lie groups/algebras.

Note that the irreducible objects in categories $\Opm$ are
precisely those irreducible $\hgl$-modules which appear
in the decomposition of $ \bFn \/(N \in \Bbb N).$ By means
of the seesaw pair above, we obtain an explicit 
formula for a decomposition of a tensor
product of two irreducible modules in
category $\Op$ (resp. $\Om$) into an infinite sum of
irreducibles in $\Op$ (resp. $\Om$) with finite multiplicities.
We further show that categories $\Opm$ are indeed semisimple
abelian tensor categories and they are equivalent to each other.
Combining with Frobenius reciprocity, we prove that
the tensor category $(\Of, \bigodot)$ is equivalent
to the tensor categories $\Opm$.

Following \cite{KR1}, we have an injective homomorphism
of Lie algebras $\widehat{\phi}_s,\, s \in \Bbb C $ 
from $\hD$ to $\hgl$
such that the pull-back of an irreducible quasifinite
$\hgl$-module via $\widehat{\phi}_s $ remains to be irreducible
as a $\hD$-module. We denote by $L_N (\hD; \kappa (\lambda), s )$
the $\hD$-module which is the pull-back via 
$\widehat{\phi}_s $ of the $N$-primitive module 
$L_N (\hgl, \Lambda_{+} (\lambda ) ), N \in \Bbb N$.
A $\hD$-module of the form of a finite tensor product
  $\bigotimes_i L_{n_i} (\hD; \kappa (\lambda^i), s_i )$
is irreducible \cite{KR1} if 
$\sum_i n_i = N, n_i \in \Bbb N,$
  $s_i \neq s_j \mbox{ mod }\Z,$ $i \neq j,$ and
  $\lambda^i \in P_{+}^{n_i}.$
We call an irreducible $\hD$-module of such form primitive.
Our category $\OWp$ of $\hD$-module  
contains all the primitive $\hD$-modules as irreducible
objects. We present an explicit formula for a tensor
product decomposition of modules in $\OWp$ and
show that the category $\OWp$ is a semisimple
abelian category.

The plan of this paper goes as follows. In Section \ref{sect_1}
we review the Lie algebras $\hgl$ and its
quasifinite hightest weight representations.
We define the two categories $\Opm$ of $\hgl$-modules.
In Section \ref{sect_duality} we present duality theorems
of the dual pair $\pairpm$ with an 
integral central charge $\pm N$.
In Section \ref{sect_decomp} we prove our 
theorem on the tensor product decomposition
in the categories $\Opm$ and show that $\Opm$ and $\Of$ are 
equivalent semisimple abelian categories.
In section \ref{sect_W} we review the Lie algebra $\hD$ and
describe
relations between quasifinite representation theory
of $\hD$ and $\hgl$. We define a category
$\OWp$ of $\hD$-modules.
In section \ref{sect_wtensor} we prove our results on
the tensor product decomposition in the category $\OWp$
and show that the category $\OWp$ is a semisimple abelian
tensor category. 

We expect that the dual pair principle \cite{H1, H2}
will have many applications in the representation 
theory of infinite dimensional Lie algebras.
Our results can be generalized
to other dual pair between a finite dimensional Lie group
and an infinite dimensional Lie algebra.
We will treat this in a future paper. We also hope our results 
will shed some lights on the constructions
of Harish-Chandra modules of $\hgl$ and $\hD$.

\section{Categories $\Opm$ of $\hgl$-modules }
\label{sect_1}

Let us denote by $\gl$ the Lie algebra of all matrices
$(a_{ij})_{i,j \in \Z}$ with only finitely many nonzero
diagonals. Letting weight $ E_{ij} = j - i $ defines
a $\Z$--principal gradation 
$\gl = \bigoplus_{j \in \Z} gl_j$.
Denote by $\nu$ the automorphism of $\gl$ which sends
$E_{ij}$ to $E_{i+1, j+1}$ $(i, j \in \Z)$.
Denote by $\hgl = gl \bigoplus \C C$ the central extension
given by the following $2$--cocycle with values in $\C$:
$$ C(A, B) = \Tr \,\left(
                     [J, A]B \right),$$
where $J = \sum_{ j \leq 0} E_{ii}$. The $\Z$--gradation
of Lie algebra $gl$ extends to $\hgl$ by letting weight $C = 0$. 
In particular, we have the triangular decomposition
$$\hgl = \widehat{gl}_{+}  
          \oplus \widehat{gl}_{0} \oplus \widehat{gl}_{-},  $$
where 
$$\widehat{gl}_{\pm} = \bigoplus_{ j \in \Bbb N}
\widehat{gl}_{\pm j}, \quad 
\widehat{gl}_{0} = gl_0 \oplus \C C.$$
Given $c \in \Bbb C$ and
$\Lambda \in {gl}_{0}^{*}$, we let
$$ \lambda_i = \Lambda (E_{ii}), \quad i \in \Z. $$
These $\lambda_i$'s are called the labels of $\Lambda$. Let
$$ H_i = E_{ii} - E_{i+1, i+1} + \delta_{i,0} C, \quad i \in \Z, $$
and
\begin{equation}
  h_i = \Lambda( H_i ) = 
   \lambda_i - \lambda_{i+1} + \delta_{i,0} c, \quad i \in \Z. 
 \label{eq_11}
\end{equation}
Denote by $L_c(\hgl, \Lambda)$ ( or $L_c(\Lambda)$ if there
is no ambiguity) the highest weight $\hgl$--module with
highest weight $\Lambda$ and central
charge $c$. Easy to see that $L_c (\Lambda)$ is quasifinite
(namely having finite dimensional graded subspaces
according to the principal gradation of $\hgl$) if
and only if all but finitely many $h_i, i \in \Z$ are zero.

Define an automorphism $\widehat{\nu}_k \; ( k \in \Bbb N)$ by
\begin{eqnarray}
 \widehat{\nu}^k (E_{ij}) &=& E_{i+k, j+k} \quad (i \neq j) \nonumber \\
 \widehat{\nu}^k (E_{ii}) &=& E_{i+k, i+k} \quad (i >0 \mbox{ or }i \leq -k) \\
 \widehat{\nu}^k (E_{ii}) &=& E_{i+k, i+k} - C \quad ( -k < i \leq 0 ) \\
 \widehat{\nu}^k ( C    ) &=& C .
                                   \nonumber
\end{eqnarray}
namely we have 
$ \widehat{\nu}_k (E_{ij}) = E_{i+k, j+k} \quad (i \neq j), 
  \widehat{\nu}_k (H_i) = H_{ i+k}.
$
Clearly we have $\widehat{\nu}^k$ is the $k$-th composition
of $ \widehat{\nu} = \widehat{\nu}^1$. One can define the
automorphism $\widehat{\nu}^k$ of $\hgl$ for $k \in - \Bbb N$
to be the inverse of $\widehat{\nu}^{-k}$. 

Define $\Lambda_j \in gl_0^* , j \in \Z$ as follows:
\begin{equation}
\Lambda_j ( E_{ii} ) =    
   \left\{
      \everymath{\displaystyle}
      \begin{array}{lll}
        1, & \mbox{for}\quad 0 < i \leq j \\
        -1, & \mbox{for}\quad j < i \leq 0  \\
        0, & \mbox{otherwise.}
      \end{array}
    \right. \\
  \label{eq_122}
\end{equation}
Define $\hL_0 \in \hgl_0^{*}$ by 
$$ \hL_0 (C) = 1, \quad
  \hL_0 (E_{ii}) = 0 \mbox{ for all } i \in \Z,$$
and extend $\Lambda_j\/$ from $gl_0^{*}\/$ to $\hgl_0^{*}$ by letting
$\Lambda_j (C) = 0.$ Then 
$$\hL_j = \Lambda_j + \hL_0, \quad j \in \Z $$ 
are the fundamental weights,
i.e. $\hL_j ( H_i ) = \delta_{ij}.$

In the case that the central charge is a non-negative integer $N$, we
call a highest weight $\Lambda$ {\em $N$-primitive} 
and a $\hgl$-module
$L(\Lambda)$ a {\em $N\/$-primitive module} if the
$h_i$'s (cf. equation(\ref{eq_11})) satisfy
\begin{equation}
h_k \in \Bbb Z_{+} \;\; (k \in \Z), \quad \sum_k h_k = N.
\end{equation}
We denote by $\Op$ the category of 
all $\hgl_0$-diagonalizable, $\hgl_{+}$-locally finite 
$\hgl$--modules, with $N$-primitive 
$\hgl$--modules for every $N \in \Z_{+} $
as all irreducible objects
and such that any module in $\Op$ has a Jordan-Holder composition
series in terms of $N\/$-primitive $\hgl$-modules
$( N \in \Z_{+})$. Denote by $\OpN$ the subcategory of $\Op$
consisting of those representations with central charge $N$.
Here and further 
a $\frak g$-module $M$ is called ${\frak g}'$-locally
finite for a Lie subalgebra $ {\frak g }' \subset {\frak g}$
if for any $a \in M,$ the ${\frak g}'$-submodule
${\cal U} ({\frak g}')a$ is finite dimensional.

In the case that the central charge is a non-positive integer $-N$, we
call a highest weight $\Lambda$ to be {\em $(-N)$-primitive}
and $L(\Lambda)$ a {\em $(-N)\/$-primitive module} if the
$h_i$'s satisfy
\begin{equation}
  h_i \in \Z_{+} \mbox{ if } i \neq 0 \mbox{ and } \sum_i h_i = -N.
\end{equation}
\begin{equation}
  \mbox{ If } h_i \neq 0 \mbox{ and } h_j \neq 0, \mbox{ then } 
    | i-j | \leq N.
\end{equation}

Denote by $\Om$ the category of 
all $\hgl_0$-diagonalizable, $\hgl_{+}$-locally finite
$\hgl$ modules, with $(-N)$-primitive 
$\hgl$-modules for every $N \in \Z_{+} $
as all irreducible objects, and such that 
any module in $\Om$ has a Jordan-Holder composition
series in terms of $(-N)\/$-primitive $\hgl$-modules
$( N \in \Z_{+})$.

\begin{remark}
   In the category $\Op$, $N$-primitive $\hgl$-modules ($ N \in \Z_{+}$)
  are exactly irreducible unitary highest weight $\hgl$-modules
  with respect to the standard compact anti-involution in $\hgl$.
\end{remark}

\section{ Duality of $(GL(N), \hgl)$ }
\label{sect_duality}

We will first describe some Fock spaces on which the 
Lie algebras $gl(N)$ and $\hgl$ act on and 
commute with each other. For $\hgl$ with
positive (resp. negative)
integral central charge $c=N$ (resp. $c=-N$), we need N pairs
of free fermions (resp. bosonic ghosts). We will
deal with these two cases in a parallel way as follows.
Let $\bb (z)$ (resp. $\bc (z)$) represent
either $b(z)$ or $\beta(z)$ (resp. $c(z)$ or $\gamma (z)$).
Introduce
$$\bb (z) = \SUM \bb \N, \quad \bc (z) = \SUM \bc \NN.$$
We have the following commutation relations
$$ [ \bc_m, \bb_n ]_{\epsilon} 
= \bc_m \bb_n+ \epsilon \bb_n \bc_m = \delta_{m+n, 0}$$
with $\epsilon  = 1$ for fermions and $\epsilon = -1$ for
bosons. We define the Fock space $\F_{+}$ (resp. $\F_{-}$)
of the fields $b(z)$ and $c(z)$ 
(resp. $\beta (z)$ and $\gamma (z)$), 
generated by the vacuum $\vac$, satisfying
$$\bc_n \vac = \bb_{n+1} \vac = 0, \quad n \geq 0.$$
We denote by $\FFn$ (resp. $\FBn$) 
the Fock space of $N$ pairs of fermions
(resp. bosonic ghosts).

Now we take $N$ pairs of $\bb\bc$ fields, 
$\bb^p (z), \bc^p (z), p = 1, \dots, N$,
and consider the
corresponding Fock space $\bFn$.

Introduce the following generating series
\begin{equation}
  E(z,w) \equiv \sum_{i,j \in \Z} E_{ij} z^{i-1}w^{-j} 
   = \epsilon \sum_{p =1}^N :\bc^p (z) \bb^p (w):
\end{equation}
\begin{equation}
  e^{pq} (z) \equiv \SUM e^{pq}\NN =\epsilon : \bc^p(z) \bb^q (z):, 
   \quad p,q = 1, \dots, N,
\end{equation}
where the normal ordering $::$ means that the operators
annihilating $\vac$ are moved to the right and multiplied by
$- \epsilon$.

It is well known (cf. e.g. \cite{FKRW, KR2})
that the operators $E_{ij},$ $ i,j \in \Z$
form a representation in $\bFn$ of
the Lie algebra $\hgl$ with central charge $\epsilon N$;
the operators $e^{pq}(n), p,q = 1, \dots, N, n \in \Z$ 
form a representation of the affine Kac-Moody algebra
$\widehat{gl(N)}$ with central charge $ \epsilon$ \cite{FF}. 
In particular the operators $e_{pq} := e^{pq}(0), \/
p, q = 1, \cdots N,$ form
the horizonal subalgebra $ gl(N)$ in $\widehat{gl(N)}$.
One can also check directly that
\begin{equation}
   [ e_{pq}, E_{ij} ] = 0, \quad
    p,q = 1, \dots, N, \quad i,j \in \Z.
\end{equation}
               
The action of $gl(N)$ on $\bFn$ can be easily shown to
lift to an action of $GL(N)$. Let 
$ P^N_{+} = \{ \lambda = (\lambda_1, \cdots, \lambda_N) \in \Z^N
\mid \lambda_1 \geq \dots \geq \lambda_N \}.$ 
Denote by ${\frak h}_N$
the diagonal matrices of the complex Lie algebra $gl(N)$. 
We denote by ${}^N V(\lambda)$ (or $V(\lambda)$ when
there is no ambiguity) the irreducible representation
of $GL(N)$ generated by a highest weight vector 
of highest weight $\lambda$ which is fixed
by the unipotent subgroup of upper triangular 
matrices with $1$'s along diagonals.

We define a map $\Lambda_{+}: P^N_+ \longrightarrow \hgl_0^*$: 
$$\lambda = (\lambda_1, \cdots, \lambda_N) 
 \longmapsto \Lambda_{+} (\lambda)$$
to be
 \begin{equation}
  \Lambda_{+} (\lambda) = \hL_{\lambda_1} + \cdots + \hL_{\lambda_N}.
   \label{map_plus}
 \end{equation}

We define another map $\Lambda_{-}: P^N_+ \longrightarrow \hgl_0^*$:
$$\lambda = (\lambda_1, \cdots, \lambda_N) 
 \longmapsto \Lambda_{-} (\lambda)$$
to be
  \begin{equation}
   \Lambda_{-} (\lambda) 
    = (\dots, 0, \lambda_{p+1}, \dots, \lambda_N; 
                 \lambda_1, \dots, \lambda_p, 0, \dots),
    \label{map_minus}
  \end{equation}
where $ p = 0$ if all $\lambda_i < 0$, $p = N $ if all
$\lambda_i > 0$, and 
$1 \leq p < N$ is chosen such that 
$\lambda_p \geq 0 \geq \lambda_{p+1}$. Here the semicolon is 
put between the 0th and the first slots and $gl^*_0$ is
identified with ${\Bbb C}^{\infty}$ in a canonical way.

Then we have the following duality theorems.

\begin{theorem}

With respect to the dual pair $(GL(N), \hgl |_{c = N})$ the module
$\FFn$ decomposes into a direct sum of irreducible
isotypic spaces:
 \begin{equation}
   \FFn = \bigoplus_{\lambda \in P^N_+} 
          V(\lambda) \otimes L_N \left(
                                 \Lambda_{+}(\lambda)
                               \right)
 \end{equation}
where $V(\lambda)$ is the irreducible highest weight $GL(N)$-module
of highest weight $\lambda$, and 
$L_N \left(  \Lambda_{+}(\lambda)  \right)$ 
(See $(\ref{map_plus})$ for the definition 
of $ \Lambda_{+}(\lambda)$) is the irreducible 
highest weight $\hgl$-module with central charge $N$.

  \label{th_dualplus}
\end{theorem}
Theorem \ref{th_dualplus} was first proved in \cite{F1} and also proved 
in \cite{FKRW} in a different method. The following theorem is proved
in \cite{KR2}.

\begin{theorem}

With respect to the dual pair $(GL(N), \hgl |_{c = N})$ the module
$\FBn$ decomposes as follows:
 \begin{equation}
   \FBn = \bigoplus_{\lambda \in P_{+}^N} 
          V(\lambda) \otimes L_{-N} \left(
                                 \Lambda_{-}(\lambda)
                               \right)
 \end{equation}
where $V(\lambda)$ is the irreducible highest weight $GL(N)$-module
of highest weight $\lambda$ and 
$L_{-N} \left(  \Lambda_{-}(\lambda)  \right)$ is the irreducible 
highest weight $\hgl$-module with central charge $-N$
(see (\ref{map_minus}) for the definition of $ \Lambda_{-}(\lambda)$).
  \label{th_dualminus}
\end{theorem}

\section{Tensor product decomposition in $\Opm$}
\label{sect_decomp}

For a given $\lambda \in P^{M+N}_{+}\/ (M, N \in \Bbb N)$, 
when we restrict the
irreducible representation ${}^{M+N} V(\lambda)$ of $GL(M+N)$
to the block-diagonal subgroup 
  \begin{eqnarray*}
  \left[ \begin{array}{cc}
    GL(M)      & 0      \\
      0        & GL(N)
         \end{array} \right]
  \end{eqnarray*}
which is identified with $GL(M) \times GL(N)$,
we have the following decomposition into a direct sum
of irreducible $GL(M) \times GL(N)$-modules
\begin{equation}
  {}^{M+N}V(\lambda) = \bigoplus_{\mu \in  P^{M}_{+}, \nu \in P^{N}_{+}}
   C^{\lambda}_{\mu \nu} {}^{M}V(\mu) \bigotimes {}^{N}V(\nu)
\label{eq_recip}
\end{equation}
where $C^{\lambda}_{\mu \nu}$ denotes the multiplicity in 
${}^{M+N}V(\lambda)$ of
the $GL(M) \times GL(N)$ module
$ {}^{M}V(\mu)  \bigotimes {}^{N}V(\nu)$ ($\bigotimes$ here
denotes the outer tensor product).

Now we are ready to state our first main result.

\begin{theorem}
The tensor product of two irreducible
$\hgl$-modules in category $\Op$ (resp. $\Om$)
can be decomposed into a direct sum of irreducible
$\hgl$-modules with
finite multiplicites. More precisely,
given $\mu \in  P^{M}_{+}, \nu \in P^{N}_{+}$, and consider
the irreducible $\hgl$-modules $L_M( \Lambda_{+}(\mu) )$ and
$L_N( \Lambda_{+}(\nu) )$, 
(resp. $L_{-M}( \Lambda_{-}(\mu) )$ and
$L_{-N}( \Lambda_{-}(\nu) )$) with central charges $M$ and $ N$
(resp. $-M$ and $-N$) in category $\Op$ (resp. $\Om$)
(cf. equations (\ref{map_plus}), (\ref{map_minus})). 
Then we have the following decomposition
  \begin{equation}
     L_{\pm M}( \Lambda_{\pm}(\mu) ) \bigotimes 
       L_{\pm N}( \Lambda_{\pm}(\nu) )
      =  \bigoplus_{\lambda \in P^{M+N}_{+} }
          C^{\lambda}_{\mu, \nu} L_{\pm (M+N)} \left(           
                                     \Lambda_{\pm}(\lambda)
                                   \right)
  \end{equation}
where the multiplicity
$C^{\lambda}_{\mu \nu}$ is defined as in (\ref{eq_recip}),
$L_{\pm(M+N)} \left( \Lambda_{\pm}(\lambda)  \right)$ 
is the irreducible $\hgl$-module in category $\Opm$
with highest weight
$\Lambda_{\pm}(\lambda)$ 
and central charge $\/\pm (M+N)$.
   \label{th_rec}
\end{theorem}
\begin{demo}{Proof}
By the duality Theorems \ref{th_dualplus}, \ref{th_dualminus},
we have the following decompositions:
 \begin{equation}
   {\cal F}_{\pm}^{\otimes (M+N)}
        = \bigoplus_{\lambda \in P_{+}^{M+N}} 
          {}^{M+N}V(\lambda) \bigotimes L_{\pm (M+N)} \left(
                                 \Lambda_{\pm}(\lambda)
                               \right),
   \label{dual_mn}
 \end{equation}
 \begin{equation}
   {\cal F}_{\pm}^{\otimes M}
        = \bigoplus_{\mu \in P_{+}^{M}} 
          {}^{M}V(\mu) \bigotimes L_{\pm M} \left(
                                 \Lambda_{\pm}(\mu)
                               \right),
   \label{dual_m}
 \end{equation}
 \begin{equation}
   {\cal F}_{\pm}^{\otimes N}
        = \bigoplus_{\nu \in P_{+}^N} 
          {}^{N}V(\nu) \bigotimes L_{\pm N} \left(
                                 \Lambda_{\pm}(\nu)
                               \right).
   \label{dual_n}
 \end{equation}

Obviously 
\begin{equation}
{\cal F}_{\pm}^{\otimes (M+N)} = {\cal F}_{\pm}^{\otimes M}
\bigotimes {\cal F}_{\pm}^{\otimes N}.
  \label{eq_same}
\end{equation}
Due to (\ref{dual_m}) and (\ref{dual_n}), in addition
to the dual pair 
\begin{equation}
  \left(
    GL(M+N), \hgl |_{c = M+N} 
  \right),
 \label{oldpair}
\end{equation}
we have another dual pair 
\begin{eqnarray}
 \left(
   GL(M) \times GL(N), 
    \hgl |_{c = M} \oplus \hgl |_{c = N}
 \right)
   \label{newpair}
\end{eqnarray}
acting on the same Fock space ${\cal F}_{\pm}^{\otimes (M+N)}$.
Clearly, we have the following inclusion relations of Lie algebras:
$$
  GL(M) \times GL(N) \subset GL(M+N), \quad
  \hgl |_{c = M+N} \subset_{diagonal}
    \hgl |_{c = M} \oplus \hgl |_{c = N}.
$$
So the two dual pairs (\ref{oldpair}) and
(\ref{newpair}) form a seesaw pair in the sense of \cite{Ku}.

>From (\ref{eq_recip}) and (\ref{dual_mn}), we have 
the following decomposition as modules over
$GL(M) \times GL(N) \times \hgl |_{ c = \pm (M+N) }$:
 \begin{equation}
      \begin{array}{rcl}
      \lefteqn{
       {\cal F}_{\pm}^{\otimes (M+N)} =
              }                       \\
      & & \bigoplus_{\lambda \in P_{+}^{M+N}} 
           \bigoplus_{\mu \in  P^{M}_{+}, \nu \in P^{N}_{+}}
            C^{\lambda}_{\mu \nu} {}^{M}V(\mu) \bigotimes {}^{N}V(\nu)
              \bigotimes L_{\pm(M+N)} \left(
                          \Lambda_{\pm}(\lambda)
                                   \right).
      \label{dual_mix1}
      \end{array}
  \end{equation}

>From (\ref{dual_m}), (\ref{dual_n}) and (\ref{eq_same}), we have
the following decomposition of the Fock space $\bFmn$
as modules over
$GL(M) \times GL(N) \times \hgl |_{ c = \pm (M+N) }$:
 \begin{equation}
      \begin{array}{rcl}
      \lefteqn{
        {\cal F}_{\pm}^{\otimes (M+N)} =       
              }                             \\
       && \bigoplus_{\mu \in P_{+}^{M}, \nu \in P_{+}^N} 
           {}^{M}V(\mu) \bigotimes {}^{N}V(\nu) 
                     \bigotimes L_{\pm M} \left(
                                 \Lambda_{\pm}(\mu)
                               \right)
                     \bigotimes L_{\pm N} \left(
                                 \Lambda_{\pm}(\nu)
                               \right).
   \label{dual_mix2}
      \end{array}
 \end{equation}

Now the theorem follows by 
comparing (\ref{dual_mix1}) with (\ref{dual_mix2}).
\end{demo}

\begin{remark}
  The multiplicities $C^{\lambda}_{\mu \nu}$ can be computed
  by the combinatorial recipe known as Littlewood-Richardson
  Rule \cite{M}. To some extend, one may view the dual pair
  $(GL(N), \hgl)$ as some kind of $ K \rightarrow \infty$
  limit of dual pairs $(GL(N), gl(K))$ with a semi-infinite twist
  and Theorem \ref{th_rec}
  is a consequence of the stability for tensor product \cite{H2},
  although the Lie algebra $\hgl$ is quite different from
  the direct limit of $ gl(K)$ as $ K \rightarrow \infty$.
\end{remark}

Note that a short exact sequence
of $\hgl$-modules
$0 \rightarrow V_1 \rightarrow V \rightarrow V_2 \rightarrow 0$,
always splits if 
$V_1$, $V_2$ are two irreducible $\hgl$ modules with different
cental charges.
Since all irreducible modules from category $\Op$ are
integrable modules \cite{K}, the complete reducibility
holds in $\Op$. Now the complete reducibility in category $\Om$
follows from Theorem 4.1 in \cite{KR2}.
Combining with Theorem \ref{th_rec}, we have
\begin{theorem}
  The categories $\Op$ and $\Om$ are semisimple abelian tensor
 categories. Moreover, we have an equivalence of the two
 tensor categories $\Op$ and $\Om$ by letting
 the irreducible object 
 $ L_N (\Lambda_{+}(\lambda)) \in \Op$
 correspond to the irreducible object
 $ L_N (\Lambda_{-}(\lambda) ) \in \Om$
 for any $N \geq 0$ and $\lambda \in P_{+}^N.$
   \label{th_semisimple}
\end{theorem}
 
Denote by $\Ofn$ the category of $GL(N)$-modules
which is a (possibly infinite) sum of 
finite dimensional irreducible modules. 
Denote by $\Of$ the direct sum of the 
categories $\Ofn$ for all $N \geq 0$, namely
a category whose objects consist of
a direct sum of modules 
of $GL(N)$ in $\Ofn$ for all $N \geq 0$.

We introduce the following tensor product $\bigodot$ on the category $\Of$:
given a module $U \in \Ofm$ and a module $V \in \Ofn$,
let $U \bigotimes V$ be the outer tensor product of $U$ and $V$
which is a $GL(M) \times GL(N)$-module. We define 
$$
  U \bigodot  V
   = \left( ind_{GL(M) \times GL(N)}^{ GL(M+N)} U \bigotimes V
     \right)^{l.f.}
$$
where $ind_H^G W$ denotes the induced $G$-module from a module $W$ of
$H \subset G$ and $X^{l.f.}$ denotes 
taking the locally finite vectors in $X$ which are by definition 
the vectors which lie in some finite dimensional $GL(M+N)$-submodule of $X$.
\begin{theorem}
  We have an equivalence of tensor categories between
 $(\Of, \bigodot)$ and $( \Opm, \bigotimes)$ by sending
 $ {}^N V(\lambda)$ to $ L_N \left(\hgl,  \Lambda_{\pm}(\lambda) \right) $.
  \label{th_equiv}
\end{theorem}
\begin{demo}{Proof}
 It suffices to prove for $\Op$ by Theorem \ref{th_semisimple}.
 It is clear that the irreducible objects in the
 two categories $\Of$ and $\Opm$ are in one-to-one correspondence.

  Take a finite dimensional irreducible $GL(M+N)$-module
 ${}^{M+N} V(\lambda)$. By Frobenius Reciprocity we have
 \begin{eqnarray*}
  \begin{array}{rcl}
   \lefteqn{
  \mbox{Hom}_{GL(M +N)}
       \left( {}^{M+N} V(\lambda), {}^M V(\mu) \bigodot  {}^N V(\nu)
       \right)      
           }              \\
  & = & \mbox{Hom}_{GL(M +N)}
      \left( {}^{M+N} V(\lambda), ind_{GL(M) \times GL(N)}^{ GL(M+N)} 
                 {}^M V(\mu) \bigotimes {}^N V(\nu)
              \right)     \\
  & = & \mbox{Hom}_{GL(M) \times GL(N)}
      \left(  {}^{M+N} V(\lambda) \mid_{GL(M) \times GL(N)},
           {}^M V(\mu) \bigotimes {}^N V(\nu)
      \right).
  \end{array}
 \end{eqnarray*}
 So we have
 $$ \mbox{dim } \mbox{Hom}_{GL(M +N)}
       \left( {}^{M+N} V(\lambda), {}^M V(\mu) \bigodot  {}^N V(\nu)
       \right)      
   = C^{\lambda}_{\mu \nu}.
 $$
 Hence we have
 $$ {}^M V(\mu) \bigodot  {}^N V(\nu)
   =  \bigoplus_{\lambda \in P^{M+N}_{+} }
          C^{\lambda}_{\mu, \nu} {}^{M+N} V(\lambda).
 $$
 The tensor category $\Of$ is clearly semisimple. Now the 
 theorem follows by comparing with Theorem~\ref{th_rec}.
\end{demo}
\begin{remark}
  It follows that the category $\Of$ is abelian as well. 
 It is not difficult to show by using Frobenius reciprocity repeatedly
 that the tensor product of $n$ ($n \geq 2$) irreducible modules 
 in the category $\Of$ (and also in $\Opm$ by Theorem \ref{th_equiv})
 can be decomposed into an (infinite) direct sum of irreducibles
 with finite multiplicities. Indeed we can define 
 a tensor subcategory of $\Of$ (resp. $\Opm$) to be the smallest
 among subcategories of $\Of$ (resp. $\Opm$) consisting of the same
 irreducible objects as in $\Of$ (resp. $\Opm$) and 
 closed under the finite direct sum and tensor product operations.
 Although an object in such a tensor subcategory may still be an infinite
 sum of irreducibles, the ``infinity'' behaves in a controlled way.
   \label{rem_sense}
\end{remark}
\begin{remark}
   Consider the automorphism $ \widehat{\nu}^{-k}$ on $\hgl$.
 Via pull-back the irreducible module
 $L( \sum_{ a =1}^N \hL_{\lambda_a})$
 in the category $\OpN$ becomes $L( \sum_{ a =1}^N \hL_{\lambda_a +k})$.
 The corresponding operation on $\Ofn$ by the equivalence of
 categories as in Theorem~\ref{th_equiv}
 is given by tensoring with $det^{ \otimes k}$, where $det$ is 
 the one-dimensional representation
 of $GL(N)$ of highest weight $(1, \ldots, 1)$. 
\end{remark}

\section{Category $\OWp$ of $\hD$-modules}
\label{sect_W}

Let $\D$ be the Lie algebra of regular differential operators on
the circle. The elements 
\begin{eqnarray*}
  J^l_k = - t^{l+k} ( \partial_t )^l, \quad
    l \in \Z_{+}, k \in \Z   \nonumber
\end{eqnarray*} 
form a basis of $\D$. $\D$ has also another basis
\begin{eqnarray*}
  L^l_k = - t^{k} D^l, l \in \Z_{+}, k \in \Z    \nonumber
\end{eqnarray*}
where $D = t \partial_t$. Denote by $\hD$ the central extension of  
$ \D $
by a one-dimensional center with a generator $C$, with
commutation relations (cf. \cite{KR1})
\begin{eqnarray}
  \left[
     t^r f(D), t^s g(D)
  \right]
    & = & t^{r+s} 
    \left(
      f(D + s) g(D) - f(D) g(D+r) 
    \right) \nonumber \\ 
   & + & \Psi 
      \left(
        t^r f(D), t^s g(D)
      \right)
      C
  \label{eq_12}
\end{eqnarray}

where 
\begin{equation}
  \Psi 
      \left(
        t^r f(D), t^s g(D)
      \right)
   = 
   \left\{
      \everymath{\displaystyle}
      \begin{array}{ll}
        \sum_{-r \leq j \leq -1} f(j) g(j+r),& r= -s \geq 0  \\
        0, & r + s \neq 0. 
      \end{array}
    \right. \\
  \label{eq_13}
\end{equation}
$\hD$ is often refered to as the $\W$ algebra in literature.

Letting weight $J^l_k = k$ and weight $ C = 0$ defines a principal
gradation
\begin{equation}
  \hD = \bigoplus_{j \in \Z} \hat{\cal D}_j
  \label{eq_14}
\end{equation}
and so we have a triangular decomposition
\begin{equation}
  \hD = \hD_{+} \bigoplus \hD_{0}  \bigoplus \hD_{-}
  \label{eq_15}
\end{equation}
where 
\begin{eqnarray*}
  \hD_{\pm} = \bigoplus_{j \in \pm \Bbb N} \hat{\cal D}_j,
  \quad
  \hD_{0} = {\D}_{0}  \bigoplus {\Bbb C} C.
\end{eqnarray*}

Given $s \in \Bbb C$, $\D$ acts naturally on the vector space 
$ t^s \Bbb C [ t, t^{-1} ] $. Choose a basis 
$ v_j = t^{s + j}, j \in \Z\,$ of the space
$\/ t^s \Bbb C [ t, t^{-1} ] $.
Then this gives rise to a monomorphism of Lie algebras 
$$ \phi_s : {{\D} {\longrightarrow} {\gl} }, \; s \in \Bbb C $$
defined by
\begin{equation}
\phi_s \left(
       t^k f(D)
     \right)
       = \sum_{ j \in \Z} f(s-j) E_{j-k,j},
  \label{eq_embed}
\end{equation}
and then a monomorphism of Lie algebras \cite{KR1}
\begin{equation}
  \widehat{\phi}_s : \hD \longrightarrow \hgl, \; s \in \Bbb C
  \label{map_phihat}
\end{equation}
defined by
\begin{eqnarray*}
  \widehat{\phi}_s \mid_{\hD_j} & = & \phi_s \mid_{ {\cal D}_j}
  \quad \mbox{if}\quad j \neq 0, \\
   \widehat{\phi}_s \left(
                     e^{xD}
                    \right)
       & = & {\phi}_s \left(
                        e^{xD}
                      \right)
              - \frac{ e^{sx} -1 }{ e^x -1 } C, \\
     \widehat{\phi} (C) & = & C.
\end{eqnarray*}

Given a sequence of complex numbers 
$\xi = (\xi_j)_{ j \in \Z_{+} }$ and a complex
number $c$, there exists a unique irreducible highest
weight $\hD$-module
$L_c (\hD, \xi) $, which admits a nonzero vector $v_{\xi}$
such that:
\begin{eqnarray*}
  L^j_k v_{\xi} = 0 \quad \mbox{for} \quad
   k > 0, \quad L^j_0 v_{\xi} = \xi_j v_{\xi},
    \quad C = cI.
\end{eqnarray*}
The module is called {\em quasifinite} if all eigenspaces
of the operator $D$ is finite dimensional.
It was proved in \cite{KR1} that $L_c (\hD, \xi)$
is a quasifinite module if and only if
the generating function
$$ \Delta_{\xi} (x) = 
     \sum_{n = 0}^{\infty} \frac{x^n}{n!} {\xi}_n $$
has the form
$$ \Delta_{\xi} (x) = 
     \frac{\phi (x) }{e^x - 1}  $$
where 
\begin{eqnarray*}
\phi (x) + c = \sum_i p_i (x) e^{r_i x} \quad (\mbox{a finite sum})
\end{eqnarray*}
for some nonzero polynomials in $x$ such that
$\sum_i p_i (0) = c$ and distinct
complex numbers $r_i$. The numbers $r_i$ are called the 
{\em exponents} of this module while 
the polynomial $p_i (x)$ are called their
{\em multiplicities}. We call an irreducible quasifinite 
$\hD$-module $L_c (\hD, \xi)$ {\em primitive} if
all multiplicities of its exponents are non-negative
integers \cite{FKRW}.
Note that a primitive $\hD$-module always has a
non-negative integral central charge. 

The following two propositions are taken from \cite{KR1}.
\begin{proposition}
We have the isomorphism
$$ L_c (\hD, \xi) \bigotimes L_{c'} (\hD, {\xi}' )
    \equiv  L_{c + c'} (\D, \xi + {\xi}' )  $$
provided that no exponent of the
first module is congruent  $\mbox{mod }\Z$ to any 
exponent of the second module.
   \label{prop_kr1}
\end{proposition}
\begin{proposition}
Let $L_c (\hgl, \Lambda)$ be the irreducible highest weight
$\hgl$-module with highest weight 
$\Lambda = ({\lambda}_j )_{ j \in \Z } $ and 
central charge $c$. $L_c (\hgl, \Lambda)$, when
regarded as a $\hD$-module via
the embedding of Lie algebras
$\widehat{\phi}_s$, is an irreducible
quasifinite $\hD$-module with exponents $s-j (\/j \in \Z)$
of multiplicity $h_j$, where $h_j$ is defined as
in $(\ref{eq_11})$.
   \label{prop_kr2}
\end{proposition}

Given $N > 0$ and $\lambda \in P^N_{+}$, $\Lp$ is the highest
weight $\hgl$-module with highest weight 
$\Lambda_{+} (\lambda ) \;(cf.\;\; (\ref{map_plus})$
for $\/\Lambda_{+} )$
and central charge $N$. We will denote the $\hD$-module
$ L_N \left(\hgl, \Lambda_{+} (\lambda ) \right)$
in the sense of Proposition \ref{prop_kr2}
by $ L_N (\hD; \kappa (\lambda), s )$
to emphasize its $\hD$-module structure. All irreducible
  quasifinite $\hD$-modules which we are concerned about
  in this paper are primitive and so have non-negative
  integral central charges. 

\begin{remark} 
   1). Any primitive $\hD$-module with
  central charge $N$ $(N \in \Z_{+})$ can 
  be realized as the form of a finite tensor product
  $$ \bigotimes_i L_{n_i} (\hD; \kappa (\lambda^i), s_i ), $$
  where $\sum_i n_i = N, n_i \in \Bbb N;$
  $s_i \neq s_j \mbox{ mod }\Z,$ $i \neq j;$ 
  $\lambda^i \in P_{+}^{n_i}.$

   2). As $\hD$-modules,
  $L_N \left(\hD; \kappa 
  (\lambda + l{\bf 1}), s + l \right), l \in \Z$
  are all isomorphic to each other, where
  ${\bf 1} = (1, \cdots, 1) \in P_{+}^N$. 
  These are all possible
  $\hD$-module isomorphisms among
  $L_N (\hD; \kappa (\lambda ), s)$, 
  $\lambda \in P_{+}^N,$ $s \in \Bbb C$ (cf. \cite{KR1}).
   \label{rem_unique}
\end{remark}

We denote by $\OWp$ the category of $\hD_0$-diagonalizable,
$\hD_{+}$-locally finite quasifinite $\hD$-modules, 
with primitive modules as all irreducible objects
and such that any module in $\OWp$ has a Jordan-Holder composition
series in terms of primitive modules.
Modules in category $\OWp$ always have
non-negative integer central charges.

\begin{remark}
  The category $\OWp$ consists of exactly all
 representations of the vertex algebra $\W$
 with non-negative integral central charges, cf. $\cite{FKRW}$.
  \label{rem_big}
\end{remark}
\begin{remark}
   If we try to study a similarly defined category
  of $\hD$-modules with non-positive integral central
  charge, we will see the complete reducibility fails.
  Indeed such a category will contain irreducible
  quasifinite $\hD$-modules 
  $L_{-1} ( \hD; \kappa(0), s ) \;\;(s \in \Bbb C)$ 
  with generating
  functions $ \Delta^s (x) = - \frac{e^{sx} - 1}{e^x - 1} $.
  But there exists  non-split short exact sequences 
  (cf. Proposition 4.2 in \cite{W}):
  $$ 0 \longrightarrow L_{-1} ( \hD; \kappa(0), s ) 
       \longrightarrow M_{\pm} 
       \longrightarrow L_{-1} ( \hD; \kappa(0), s \pm 1 ) 
       \longrightarrow 0 $$ 
  for some $\hD$-module $M_{ \pm}$.
  Equivalently, there
  exists non-split short exact sequences of
  $\hgl$-modules with central charge $-1$:
  $$ 0 \longrightarrow L_{-1}(\hgl, - \widehat{\Lambda}_i)
      \longrightarrow V_{\pm } 
      \longrightarrow L_{-1}(\hgl, - \widehat{\Lambda}_{i \pm 1} ) 
      \longrightarrow 0, $$ 
 for some $\hgl$-module $V_{\pm }$,
  where $\widehat{\Lambda}_j$ means the $j$-th fundamental weight
  of $\hgl$.
\end{remark}

\section{Tensor product decomposition in category $\OWp$}
 \label{sect_wtensor}

Take two irreducible $\hD$-modules $V_1$ and $V_2$ in
the category $\OWp$.
First let us consider the simplest case when 
all exponents in $V_1$ (resp. $V_2$)
belong to the same congruence class $\mbox{mod } \Z$.
Then by Remark \ref{rem_unique}, we may express
$V_1$ and $V_2$ in the following forms:
$$ V_1 = L_{ M}(\hD; \kappa  (\mu), s ),
    \quad V_2 = L_{N}(\hD; \kappa  (\nu), t ).$$
Now we have two possibilities.

First, if the congruence class of exponents 
in $L_M(\hD; \kappa  (\mu), s )$ does
not coincide with the one in $L_{N}(\hD; \kappa  (\nu), t )$,
then the tensor product
$$L_{M}(\hD; \kappa  (\mu), s ) \bigotimes 
       L_{N}(\hD; \kappa  (\nu), t )        $$
is an irreducible $\hD$-module in category $\OWp$
by Proposition \ref{prop_kr1}.

Secondly, if the congruence class of exponents in 
$L_{N}(\hD; \kappa  (\mu), s )$
coincides with the one in $L_{N}(\hD; \kappa  (\nu), t )$, 
we may assume $s = t$ by Remark \ref{rem_unique}. Then by 
Theorem \ref{th_rec} and Proposition \ref{prop_kr2},
we have the following tensor product decomposition 
$V_1 \bigotimes V_2$ into a direct sum of irreducible
$\hD$-modules in category $\OWp$ 
(see (\ref{eq_recip}) for notations):
\begin{equation}
  L_{M}(\hD; \kappa  (\mu), s ) \bigotimes
   L_{N}(\hD; \kappa  (\nu), s ) = 
    \bigoplus_{ \lambda \in P^{M + N}_{+} }
       C_{ \mu \nu}^{\lambda} L_{M+N}(\hD; \kappa  (\lambda), s )
   \label{eq_basic}
\end{equation}

In general, let $s_1, s_2, \cdots $ 
(resp. $t_1, t_2, \cdots $) be representatives
of different $\mbox{mod } \Z$ congruence classes of the module
$V_1$ (resp. $V_2$) with non-negative integral central charge 
$M$ (resp. $N$). 
By Remark \ref{rem_unique}, we may write
$V_1$ and $V_2$ in the following forms:
\begin{equation}
V_1  =
 \bigotimes_i L_{m_i} (\hD; \kappa (\mu^i), s_i ),
  \label{eq_v1}
\end{equation}
where $\sum_i m_i = M, m_i \in \Z_{+};
\mu^i \in P_{+}^{m_i};$
and
\begin{equation}
V_2 =
 \bigotimes_j L_{n_j} (\hD; \kappa (\nu^j), t_j ),
  \label{eq_v2}
\end{equation} 
where $\sum_j n_j = N, n_j \in \Z_{+};
\lambda^j \in P_{+}^{n_j}.$
We may assume (by rearrangement of the order
of $s_i$'s and $t_j$'s if necessary) $s_a = t_a \mbox{ mod } \Z,$ 
$a = 1, \cdots, k $ for some $ k \geq 0.$ 
Here $k = 0$ would mean that none of the exponents $s_1, s_2, \cdots $ 
is congruent $\mbox {mod }\Z$
to any of the exponents 
$t_1, t_2, \cdots $. In this case, the tensor product
$V_1 \bigotimes V_2$
is an irreducible $\hD$-module in $\Op$ 
by Proposition \ref{prop_kr1}.

We may further assume $t_a = s_a, a = 1, \cdots, k$ 
by Remark \ref{rem_unique}. Then we have the 
following tensor product decomposition 
into a direct sum of irreducible $\hD$-modules
in $\OWp$ according to (\ref{eq_basic}):
\begin{eqnarray}
  L_{m_a}(\hD; \kappa  (\mu^a), s_a ) \bigotimes
   L_{n_a}(\hD; \kappa (\nu^a), s_a )   \nonumber \\
   =
    \bigoplus_{ \lambda^a \in P^{m_a + n_a}_{+} }
       C_{ \mu^a \nu^a}^{\lambda^a} L(\hD; \kappa  (\lambda^a), s_a),
  \quad a = 1, \cdots, k.
   \label{eq_single}
\end{eqnarray}
                    
Now we are ready to state our main result of this section.
\begin{theorem}
With notations as in $(\ref{eq_v1}, \ref{eq_v2}, \ref{eq_single})$, 
we have the following decomposition
of a tensor product of two modules $V_1$ and $V_2$
in $\OWp$ into the direct sum of irreducibles in $\OWp$:
\begin{eqnarray}
 V_1 \bigotimes V_2 
  & = &  \bigoplus_{ \lambda^1 \in P^{m_1 + n_1}_{+} }
       \cdots
       \bigoplus_{ \lambda^k \in P^{m_k + n_k}_{+} }
       \{ \Pi_{a=1}^k C_{ \mu^a \nu^a}^{\lambda^a} \}\cdot \nonumber  \\
 & &   \cdot \left\{
        \bigotimes_{b=1}^k L(\hD; \kappa  (\lambda^b), s_b)
        \bigotimes
        \bigotimes_{i > k} L_{m_i} (\hD; \kappa (\mu^i), s_i ) 
                                                \right.\nonumber  \\
 & &    \left. \bigotimes
        \bigotimes_{j > k} L_{n_j} (\hD; \kappa (\nu^j), t_j )
       \right\}.
       \label{eq_messy}
\end{eqnarray}
\end{theorem}
\begin{demo}{Sketch of a proof}
  By (\ref{eq_single}), Propositions \ref{prop_kr1}, \ref{prop_kr2} and
Remark \ref{rem_unique}, we have
\begin{eqnarray}
 V_1 \bigotimes V_2 
 &=& 
     \bigotimes_{1 \leq a \leq k} 
       \bigoplus_{ \lambda^a \in P^{m_a + n_a}_{+} }
        \left \{ 
           C_{ \mu^a \nu^a}^{\lambda^a} 
            L(\hD; \kappa  (\lambda^a), s_a)
        \right \} 
     \bigotimes \nonumber  \\
 & &  
     \bigotimes_{i > k} L_{m_i} (\hD; \kappa (\mu^i), s_i )
     \bigotimes
     \bigotimes_{j > k} L_{n_j} (\hD;\kappa (\nu^j), t_j ).
       \label{eq_plug}
\end{eqnarray}
Then it is easy to see that the right hand side of 
(\ref{eq_plug}) is equal to the right hand side of
(\ref{eq_messy}). Each of the $\hD$-modules appearing
on the right hand side of (\ref{eq_messy}) is
irreducible by Proposition \ref{prop_kr1}
and is in category $\OWp$ by Remark \ref{rem_unique}.
\end{demo}

Complete reducibility in Category $\OWp$
follows from complete reducibility in Category $\OWp$
with the help of Propositions \ref{prop_kr1},
\ref{prop_kr2}, and Remark \ref{rem_unique}. So we 
have proved the following theorem.

\begin{theorem}
   $\OWp$ is a semisimple abelian tensor category.
\end{theorem}
\begin{remark}
 A tensor subcategory of $\OWp$ can be defined analogous to
 the tensor subcategory of $\Opm$ defined in Remark \ref{rem_sense}.
\end{remark}

{\bf Acknowledgement} I thank Igor Frenkel,
Roger Howe and Yan Soibelman for stimulating discussions.

\noindent
E-mail address: {\tt wqwang@math.yale.edu}

\end{document}